\documentclass[11pt]{article}
\usepackage{latexsym}
\usepackage{amsfonts}
\usepackage{graphicx}
\usepackage{epstopdf}
\usepackage{color}
\usepackage[sort]{natbib}
\usepackage{lineno}

\newtheorem{theorem}{Theorem}
\newtheorem{lemma}{Lemma}
\newtheorem{corol}{Corollary}
\newtheorem{property}{Property}

\newcommand{\be}{\begin{equation}}
\newcommand{\ee}{\end{equation}}

\newcommand{\narrow}{\setlength{\itemsep}{1pt}
	\setlength{\parskip}{0pt}\setlength{\parsep}{0pt}}
\begin{document}
	\title{\bf Improved Starting Solutions for the Planar $p$-Median Problem}
	\author{Jack Brimberg\\Department of Mathematics and Computer Science\\The Royal Military College of Canada\\ Kingston, ON Canada.\\ e-mail: Jack.Brimberg@rmc.ca\and Zvi Drezner\\
	Steven G. Mihaylo College of Business and Economics\\
	California State University-Fullerton\\
	Fullerton, CA 92834.\\e-mail: zdrezner@fullerton.edu
}
	\date{}
	\maketitle
	
\begin{abstract}

In this paper we present two new approaches for finding good starting solutions to the planar $p$-median problem. Both methods rely on a discrete approximation of the continuous model that restricts the facility locations to the given set of demand points. The first method adapts the first phase of a greedy random construction algorithm proposed  for the minimum sum of squares clustering problem. The second one implements a simple descent procedure based on vertex exchange. The resulting solution is then used as a starting point in a local search heuristic that iterates between the well-known Cooper's alternating locate-allocate method and a transfer follow-up step with a new and more effective selection rule. Extensive computational experiments show that (1) using good starting solutions can significantly improve the performance of local search, and (2) using a hybrid algorithm that combines good starting solutions with a ``deep" local search can be an effective strategy for solving a diversity of planar $p$-median problem.
\end{abstract}
\noindent{\it Key Words: Facility Location; Planar p-median; Continuous Location.}
	
	\renewcommand{\baselinestretch}{1.6}
	\renewcommand{\arraystretch}{0.625}
	\large
	\normalsize
	
\section{Introduction}
	
	The continuous $p$-median problem \citep{DBMS13,DS15}, also known as the multi-source Weber problem \citep{BHMT00,KS72} or unconstrained continuous location-allocation problem, requires generating the sites of a given number ($p$) of facilities in Euclidean space in order to minimize a weighted sum of distances from a set of demand points (fixed points, existing facilities, customers) to their respective closest facilities. Let $X_i$ denote the location of facility $i \in I= \{1,\ldots, p\}$, and $A_j$ the known location of demand point $j \in J=\{1,\ldots, n\}$. In the typical scenario, which is assumed here, the $X_i$ and $A_j$ are points in the plane, such that $X_i = (x_i, y_i)$ and $A_j = (a_j, b_j)$ for all $i, j$. As well, distances are assumed to be measured by the Euclidean norm, so that
	\begin{equation}\label{Beq1}
	d (X_i, A_j) = \sqrt{(x_i - a_j)^2 + (y_i - b_j)^2}.
	\end{equation} 
	
	Letting the weights (or demands) at the $A_j$ be given by $w_j > 0, j \in J$, we formalize the problem as follows:
	\begin{equation}\label{Beq2}
\min\left\{ f (X) = \sum\limits_{j = 1}^n w_j \min\limits_{1\le i\le p}\{ d (X_i, A_j)\}\right\},
	\end{equation}
where $X = \{X_1,\ldots, X_p\}$ denotes the set of location variables.
	
This model was originally proposed by \citet{Co63,Co64}, who also observed that the objective function $f(X)$ is non-convex, and may contain several local optima. The problem was later shown to be NP-hard \citep{MegSup}.  For historical review of the planar $p$-median problem see, for example,  \cite{LMW88,BH11}.

The single facility 1-median problem termed the Weber problem \citep{W09} has a long history dating back to the French Mathematician Pierre Fermat of the 1600s. Recent reviews of the Weber problem are \citet{W93,C19,DKSW02}.	

It is well recognized that providing a good starting solution can improve the performance of heuristic solution methods. Also, good solutions to a discrete approximation of the planar $p$-median problem that uses the set of demand points as the candidate sites can yield good solutions to the continuous model \citep[e.g., see][]{Co63,Co64,BHMT00,HMT98}. With this in mind we propose two new approaches for finding good starting points. The first method adapts phase 1 of a greedy random construction algorithm developed in \cite{KBD19} for the minimum sum of squares clustering problem. A total of $p$ dispersed demand points are identified that are likely to belong to different clusters, and these points then serve as the initial facility locations. The second method implements a simple descent procedure that starts with a random selection of $p$ demand points, and proceeds to a local minimum using vertex exchange \citep{TB68}. A third approach applies the two preceding ones in sequence. The three types of starting solutions above, and a fourth with randomly selected demand points, are all improved by the same local search derived from the IALT algorithm in \citet{BD12}. We also propose a new selection rule for the transfer follow-up step in IALT that proves to be more effective on average than the existing one. Extensive computational experiments with the four developed algorithms demonstrate that using good starting solutions can significantly improve the performance of local search in the planar $p$-median problem. The three hybrid algorithms that combine good starting solutions and ``deep" local search also prove to be quite competitive with the latest meta heuristic approaches.

The contributions of the paper include: (i) a construction algorithm for finding a good starting solution, (ii) an improvement to Cooper's  exchange algorithm ALT termed RATIO that significantly improved the IALT \citep{BD12} algorithm (which is an improved ALT), (iii) two new best known solutions to an extensively researched data set. These new approaches performed very well on a set of test problems.
	
\section{Cooper's ALT Algorithm}

Of the various heuristics proposed by \citet{Co63,Co64}, one stands out as the most famous, and is often referred to in the literature as Cooper's algorithm. This method is based on the following fundamental insight: (i) if we fix the locations of the facilities, each demand point may be allocated directly to its closest facility (with ties broken arbitrarily); (ii) fixing the resulting allocations (or partition of the customer set) results in $p$ independent single facility problems, which are easy to solve because the objective function of each sub-problem is convex. Thus, the algorithm alternates between the two phases, location and allocation, all the while reducing the objective function, until a local minimum is reached. Different variations of Cooper's method have been proposed in the literature, for example by using different starting procedures. A standard approach, commonly applied in a multi-start local search and referred to here as ALT (for “Alternating”), is outlined in Algorithm 1 \citep[see, e.g., ][]{BH11}.

\subsection*{Algorithm 1 (Standard Cooper Algorithm, ALT)}
\begin{enumerate}
	\item \label{sta1} Select the starting locations of the facilities, $X^{(0)} = \{ X_1^{(0)},\ldots, X_p^{0)}\}$, where the facility locations $X_i$ are random points in the convex hull of the given (fixed) points (or the smallest rectangle parallel to the axes that contains the set of fixed points).
 Assign each fixed point to its closest facility (with ties broken arbitrarily), resulting in a partition of the customer set $S^{(0)} = \{S_1^{(0)},\ldots, S_p^{(0}\}$, where $S_i$ denotes the subset of points assigned to facility $i \in I=\{1,\ldots, p\}$.  The starting solution is given by $\{X^{(0)}, S^{(0)}\}$. Set the iteration counter $t = 0$.
	\item\label{loc} (location phase). Solve $p$ independent single facility problems given by:
	$$\min\left\{ g_i(X_i) = \sum\limits_{j\in S_i^{(t)}} w_j d(X_i, A_j)\right\} , i = 1,\ldots, p.$$
	Let the resulting median points be given by $X^{(t + 1)} =\left\{ X_1^{(t + 1)},\ldots, X_p^{(t + 1)}\right\}$.
	\item (allocation phase). Assign each fixed point to its closest facility, using new locations $X^{(t + 1)}$; let $S^{(t + 1)} = \left\{S_1^{(t + 1)},\ldots, S_p^{(t + 1)}\right\}$ denote the new partition.
	\item (detecting a local optimum). If $S^{(t + 1)} = S^{(t)}$, stop; else set $t = t + 1$ (end of current iteration) and return to step \ref{loc}.
\end{enumerate}

\subsection{Simple Modifications}
\begin{enumerate}
	\item 	Selecting random points in the convex hull or smallest rectangle enclosing the set of demand points in Step \ref{sta1} of ALT may lead to local solutions that are degenerate, in the sense that some of the facilities will have no demand points assigned to them; see \cite{BS19}. These solutions are typically of very poor quality and should be avoided. To reduce the probability of obtaining a degenerate solution we randomly select $p$ demand points as the initial facility locations in Step \ref{sta1}.
	\item To improve the efficiency of the algorithm we only solve the location problem in Step \ref{loc} for those facilities $i$ whose allocated set $S_i^{(t)}$ has changed.
\end{enumerate}

\section{\label{Ex}Examples}

To illustrate the issues which may complicate solution approaches we analyze two simple examples with equal weights.
	
	\subsection{A Single Rectangle}
	
\begin{figure}[ht!]
		\begin{center}
			\setlength{\unitlength}{1in}
			\begin{picture}(1.5,2)
			\put(0,0){\circle*{0.08}}
			\put(2,0){\circle*{0.08}}
			\put(0,2){\circle*{0.08}}
			\thicklines
			\put(0,0){\line(1,0){2}}
			\put(0,0){\line(0,1){2}}	
			\put(0,2){\line(1,-1){2}}
			\put(0,0){\line(1,1){0.6667}}
			\put(0.6667,0.6667){\line(2,-1){1.333}}
			\put(0,2){\line(1,-2){0.6667}}
			\put(0.6667,0.6667){\circle{0.08}}
			\put(0.55,0.45){$120^\circ$}
			\put(0.7,0.7){$120^\circ$}
			\put(0.3,0.6){$120^\circ$}
			\put(0.12,0.02){$\theta$}
			\put(1.25,0.02){$60^\circ-\theta$}
			\put(-0.62,0.05){$90^\circ-\theta\rightarrow$}
			\put(-0.62,1.8){$\theta-30^\circ\rightarrow$}
			\put(-0.1,0.95){1}
			\put(1,-0.15){$a$}
			\put(0.35,0.3){$x$}
			\put(1.17,0.3){$y$}
			\put(0.28,1.2){$z$}
			\end{picture}
			\caption{\label{triang}The triangle}
		\end{center}
	\end{figure}
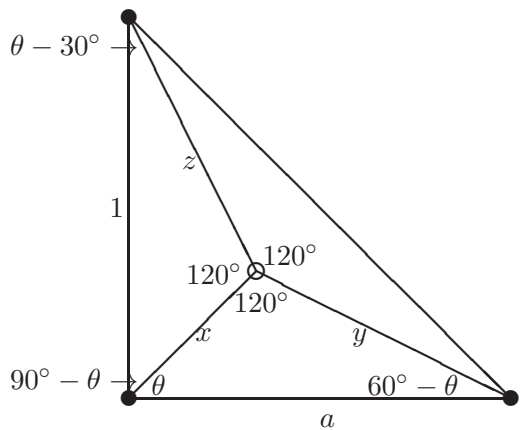

	Consider a rectangle of sides $a>1$ by 1. Demand points are located at the four corners of the rectangle ($n=4$); $w_j=1~\forall j$; and $p=2$ facilities are to be located. Locating them at the  centers of the ``longer" sides would be a terminal solution of the ALT algorithm. The objective is $2a$, while locating facilities at the shorter sides (also a terminal solution) has an objective of 2. Consider moving one point to another cluster forming clusters of 3 points and one point. The objective of the one point is 0. The objective of the triangle is more difficult to evaluate, see Figure \ref{triang}. Cavalieri proved in 1647 that the angles between the lines connecting the optimal point to the vertices of the triangle are 120$^\circ$ each \citep[see ][page 2]{DKSW02}.

	By the sinuses theorem:
	$$
	\frac{a}{\sin[120^\circ]}=\frac{x}{\sin[60^\circ-\theta]}=\frac{y}{\sin[\theta]};~~\frac{1}{\sin[120^\circ]}=\frac{x}{\sin[\theta-30^\circ]}=\frac{z}{\sin[90^\circ-\theta]}~.
	$$
	Equating the value for $x$ we get $a\sin[60^\circ-\theta]=\sin[\theta-30^\circ]$ yielding $\tan\theta=\frac{1+a\sqrt3}{\sqrt3+a}$.
	
The value of the objective function is $x+y+z$ yielding after algebraic and trigonometric manipulations a value of 
$\sqrt{a^2+a\sqrt3+1}$. We verified this result by the Solver in Excel. For a square, $a=1$, the objective in the triangle is $\sqrt{2+\sqrt3}=1.932$. It is interesting that the objective function is less than 2 for $1\le a< \frac{\sqrt{3}}{2}(\sqrt5-1)\approx 1.0705$, and thus the 3-1  clusters are optimal for the 2-median problem for these values of $a$, and not the 2-2 clusters!  Even if we start with a solution of 2 facilities at the centers of the short sides, the ALT algorithm terminates without improvement, while moving one point to the other cluster yields the optimum.
	
Consider also the following interesting point. Suppose that two vertices are selected as a starting solution for the ALT algorithm that solves the $a>1$ problem. There are six pairs that can be selected. There is a  $\frac23$ probability that the the assignment of points to clusters is on the two short sides, and $\frac13$ probability that the long sides are selected. Those are final ALT solutions. For $a>1.07$ there is a probability of $\frac13$ that a non-optimal solution is found and for $1\le a\le 1.07$ the probability of finding an optimal solution by the ALT algorithm is 0\%. However, evaluating a single point transfer will yield the optimal solution 100\% of the time.
Consider the case of a pair of points selected on the vertices of the long sides. The objective is $2a$ and creating a triangle as in Figure~\ref{triang} yields an objective of $\sqrt{a^2+a\sqrt3+1}<\sqrt{a^2+2a+1}=a+1<2a$. For $a<1.07$ this transfer leads to the optimum, and for $a>1.07$ another transfer yields the optimum.

	\subsection{Two Rectangles and Manhattan Distance}
	
\begin{figure}[ht!]
		\setlength{\unitlength}{1in}
		\centering	\begin{picture}(6,1.2)
		\put(0,0){\circle*{0.1}}
		\put(1.6,1){\circle*{0.1}}		
		\put(0,1){\circle*{0.1}}		
		\put(1.6,0){\circle*{0.1}}		
		\put(4.1,0){\circle*{0.1}}		
		\put(4.1,1){\circle*{0.1}}		
		\put(5.7,0){\circle*{0.1}}		
		\put(5.7,1){\circle*{0.1}}
		\put(1.6,1){\line(1,0){2.5}}
		\put(2.8,1.05){$x$}
		\put(0,0){\line(0,1){1}}
		\put(1.6,0){\line(0,1){1}}
		\put(4.1,0){\line(0,1){1}}
		\put(5.7,0){\line(0,1){1}}
		\put(0,0){\line(1,0){1.6}}
		\put(4.1,0){\line(1,0){1.6}}
		\put(0,1){\line(1,0){1.6}}
		\put(4.1,1){\line(1,0){1.6}}
		\put(-0.1,0.45){1}
		\put(5.75,0.45){1}
		\put(0.75,1.05){1.6}
		\put(4.8,1.05){1.6}
		\put(-0.05,-0.17){$A$}
		\put(1.55,-0.17){$B$}
		\put(4.05,-0.17){$C$}
		\put(5.65,-0.17){$D$}
		\put(-0.05,1.07){$E$}
		\put(1.55,1.07){$F$}
		\put(4.05,1.07){$G$}
		\put(5.65,1.07){$H$}
		\end{picture}
		\caption{The two rectangles example}
		\label{16}
	\end{figure}
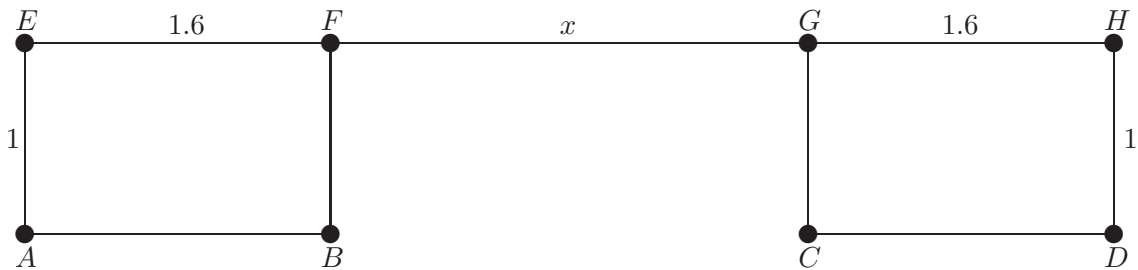

Another interesting example is given by a set of $n=8$ demand points located at the vertices of two rectangles, and the location of $p=2$ facilities using $\ell_1$ distances also called Manhattan distances \citep{LMW88,FraMcgWhi}: $d_{ij}=|x_i-x_j|+|y_i-y_j|$. The example is depicted in Figure~\ref{16} with dimension $x=2.5$.
	
The optimal location of a facility in a rectangle is anywhere in the rectangle with an objective of $2.6\times 2=5.2$, so the objective for the two rectangles is 10.4. Consider the location of the facilities at points $G$ and $A$. The objective is 10.4. However, when $F$ is transfered to the right cluster, the optimal locations of the facilities remain at $A,G$,  and the distance for $F$ is changed from 2.6 to $x$; hence the 5-3 clusters have an objective change of $x-2.6$ yielding a total of $7.8+x$, which is better for $x<2.6$. Transferring $B$ (after $F$ was transferred) to the right cluster changes the distance of the transferred point  from 1.6 to $1+x$,  a change of $x-0.6$ which further improves the objective function for $x<0.6$. Another clustering is the top 4 points as a cluster paired with the bottom 4 points. Its objective is $6.4+4x$. The optimal solution as a function of $x$ is:
	
	\bigskip
	
	\begin{tabular}{|c|c|c|}
		\hline
		Clusters&Objective&Optimal for\\
		\hline
		ABEF$|$CDGH&10.4&$x\ge 2.6$\\
		ABE$|$FCDGH&$7.8+x$&$0.6\le x\le 2.6$\\
		AE$|$BFCDGH&$7.2+2x$&$0.4\le x\le 0.6$\\
	ABCD$|$EFGH&$6.4+4x$&$0\le x\le 0.4$\\
		\hline
	\end{tabular}
	
	\bigskip
	
Visually, in Figure \ref{16} there are two clear clusters (if we think in terms of Euclidean distance) consisting of the two rectangles. However, by the analysis above transferring one of the four points $F,B,G,C$ to the other rectangle improves the value of the objective function. This shows that the $\ell_1$ measure may lead to quite different (and odd-looking) clusters.
	
For squared Euclidean distances, only the two rectangles and the 6-2 clusters can be optimal. The two rectangles objective is 7.12 and the 6-2 clusters has an objective of $\frac13 (4x^2+6.4x+16.24)$. The 6-2 clusters are optimal for $x\le \sqrt{1.92}-0.8=0.58564$ and the two rectangles are optimal for larger values of $x$.
	
\section{Other Solution Approaches}

The following algorithms are fast and therefore can be replicated many times. More elaborate algorithms applying tabu search \citep{GL97}, genetic algorithms \citep{H75,G06}, variable neighborhood search \citep{HM97, MH97}, and other meta-heuristics incorporate such fast algorithms in their approach. For a review of meta heuristic approaches to the planar $p$-median problem see \cite{DTD18,DS15,BHMT00,DBMS14,BHM08}.

\subsection{The IALT Algorithm}

The IALT modification was proposed in \cite{BD12}. The main change in the algorithm is the addition of a second phase of local search once the ALT algorithm terminates, where moves that assign a demand point from its closest facility to a farther facility are considered. The examples in Section \ref{Ex} illustrate that once the ALT algorithm terminates, improved solutions may be obtained by moving points from the closest facility to a farther facility, thus changing the set of demand points  associated with each facility and consequently the locations of the facilities. A parameter $L$ for the maximum number of transfers considered, is given. In \citet{BD12} $L=20$ was proposed.

\begin{enumerate}
\item Once the ALT algorithm terminates, for each demand point $j\in J$ the difference $\delta_j$ between the distance to the closest facility and the second closest is calculated.
\item The $L$ smallest values of $\delta_j$ are selected for transfer starting at the smallest $\delta_j$ and continuing in order.
\item The $L$ transfers to be examined are transferring the assignment of a demand point from its closest facility to its second closest.
\item The two affected facilities are relocated.
\item If the objective function does not improve, the points are returned to their original clusters and both facility  locations are restored.
\item If a transfer leads to an improved value of the objective function, it is performed, and the ALT algorithm re-started.
\item If all $L$ transfers fail to improve the value of the objective function, the algorithm terminates.
\end{enumerate}

For complete details see \cite{BD12}.

\subsection{The RATIO Modification to IALT}

Most improvement algorithms for the continuous $p$-median problem  are variations on the ALT algorithm. When each demand point is allocated to the closest facility, which in turn is optimally located with respect to its assigned cluster, the algorithm terminates. Such a final solution is a local optimum because transferring a point to another cluster without changing the facility locations  increases the value of the objective function. See also the examples in Section \ref{Ex}.
\citet{KBD19} investigated the unweighted squared Euclidean $p$-median problem. They showed the counter-intuitive property that if a point is transferred from the closest cluster of $m_1$ points to a farther cluster of $m_2$ points, the objective function {\it must decrease} if the squared distance to the second cluster center is less than $\frac{1+\frac1{m_{2}}}{1-\frac1{m_{1}}}$ times the squared distance to the closer cluster center. This ratio is, of course, different for the $p$-median problem. In the $p$-median problem  the ratio may even equal 1, meaning that the transfer cannot improve the objective function. For example, if the facility is located at a demand point, its location may not change when adding or removing a demand point from its cluster.

Consider the examples in Section \ref{Ex}. In the rectangle example, $m_1=m_2=2$, and therefore if the ratio of the squared distances is less than 3 the transfer {\it must} yield a lower squared Euclidean objective. When the rectangle is a square ($a=1$), Euclidean distances (not squared Euclidean) apply, and the facilities are located at centers of opposite sides, the distance to the closest facility is 0.5 and the distance to the second facility is $\sqrt{1.25}$ for a ratio of $\sqrt{5}\approx 2.236$. In this case we saw that the objective function is actually improved by transferring a point to the farther cluster. The ratio of the square of the distances is 5, and thus the clusters for squared distances remain at 2-2, and not 3-1, where the objective would increase from 1 to $\frac43$. For general $a\ge 1$ the objective for the two facilities at the centers of the long sides is $a^2$, and the 3-1 objective is $\frac23(a^2+1)$, which is better for $a>\sqrt2$. The ratio between the squared distances is $\frac{\frac14a^2+1}{\frac14a^2}$ which is less than 3 for $a>\sqrt2$ as shown in \cite{KBD19}.
If any value of $a$ is considered, and the facility location are in the middle of the sides of size 1, the objective is equal to 1, the squared distances to the facility on the same side are 0.25 each, and to the other side they are $0.25+a^2$ yielding a ratio of $1+4 a^2$, which is less than 3 for $a<\frac12\sqrt2$, leading to the same conclusion.

\begin{description}
\item[The RATIO Modification:]

Instead of examining all possible exchanges of a single demand point from one cluster to another (i.e., the entire one-exchange neighborhood), the IALT algorithm selects only $L$ candidate transfers based on the smallest differences  between the distances of a demand point to its second-closest and its currently assigned closest facilities. In this way the selection rule attempts to identify a short list of transfers that are most likely to improve the solution, and hence, eliminate unnecessary computations. Based on the discussion above, and the properties shown in  \cite{KBD19}, we tested the ratio between distances rather than the difference. This requires changing only one simple command in the FORTRAN code. We term IALT with this modification as the ``RATIO" variant.
\end{description}

\section{Suggestions for Generating Good Starting Solutions}

Most algorithms use randomly generated starting solutions for the facilities. Such starting solutions are usually constructed as random points in the smallest rectangle that includes all the demand points. Some algorithms use $p$ randomly selected demand points as the starting locations for the facilities. We propose simple and fast approaches for an improved selection of demand points as starting solutions.

\begin{description}
	\item[Random Selection (RAND):] This is the approach which will be used here for comparison purposes. Randomly select $p$ different demand points.
\item[The Construction Approach (CONS):]

The construction approach is similar to the first phase of the proposed construction algorithm in \citet{KGD19,KBD19}. It incorporates the idea behind the ``Greedy Randomized Adaptive Search Procedure" (GRASP) suggested by \citet{FR95}.	

\begin{enumerate}
	\item The first two demand points are randomly selected.
	\item\label{stp2} The demand point with the largest minimum distance to the already selected points is selected with probability $\frac23$ and the one with the second-largest minimum distance is selected with probability $\frac13$.
	\item Repeat Step \ref{stp2} until $p$ demand points are selected. 
\end{enumerate}

\item[The Descent Approach (DESC):]
Here we propose a minor modification of the approach in \cite{BDMS13,BD12} in order to speed up the procedure for obtaining a ``good" starting solution. That is, instead of examining the vertex swap neighborhood on a rectangular grid of $O(n^2)$ points obtained by drawing horizontal and vertical lines through the demand points as in \cite{BDMS13,BD12}, we restrict the swap to the $n$ demand points, in effect performing a local search similar to \cite{TB68},  for the ``discrete" $p$-median problem.
The second one, which is similar to the approach in \cite{DR92a}, was found to be faster.

\subsection*{Algorithm 3a}

\begin{enumerate}
	\item Randomly select $p$ demand points out of all $j\in J$ and evaluate the value of the objective function.
	\item \label{ss2}Check all $p(n-p)$ combinations of replacing a selected demand point with a non-selected one.
	\item If an improving exchange is found, select the best improving exchange and go to Step \ref{ss2}.
	\item If no improving exchange is found, stop.
\end{enumerate}

\subsection*{Algorithm 3b}

\begin{enumerate}
	\item Randomly select $p$ demand points from the set $J$ to form the set $P$,  and evaluate the objective function $f(P)$.
	\item \label{ss2}  Select one by one all demand points $i\in P$ for removal in random order, forming the set $P^\prime=P\setminus \{i\}$.
	\item \label{ss3}Check all $n-p$ points $j\in J\setminus P$ in random order for adding to the set $P^\prime$.
	\item \label{ss4} Let $P^{\prime\prime}=P^\prime \cup \{j\}$. Calculate $f(P^{\prime\prime})$.
	\item If $f(P^{\prime\prime})<f(P)$, set $P=P^{\prime\prime}$, $f(P)=f(P^{\prime\prime})$, and go to step \ref{ss2}.
	\item Otherwise,
	\begin{enumerate}
		\item select the next $j$, and go to step \ref{ss4};
		\item if all $j\in J\setminus P$ are selected, move to the next $i\in P$, set $P^\prime=P\setminus \{i\}$ and go to step~\ref{ss3};
		\item if all $i\in P$ are evaluated, stop.
	\end{enumerate}
\end{enumerate}
We can try to improve the solution further by tabu search and other heuristics but finding the best selection of $p$ points is not crucial when IALT or RATIO are applied on the selection. A ``good" one will do. Also, if the procedure is repeated many times, we would like to get different selections each time for a rich diversity of good starting solutions.

\item[The Combined Approach (COMB):]
Start the descent approach (DESC) with the construction (CONS) solution rather than a random start (RAND).

\end{description}

\subsection{\label{short}A Time Saving Short-Cut for DESC (and COMB)}

A simple and straightforward calculation of the objective function once $p$ demand points are selected for facility locations, is to go over all $n$ demand points and find the shortest distances to the $p$ selected facilities, and then sum up the weighted shortest distances. This requires $O(np)$ operations. When evaluating the objective function in a vertex exchange descent search, the complexity can be reduced to $O(n)$ \citep{Wh83} saving run time by about a factor of $p$. In our computational experiments we found that run time with the short-cut for the $p=25$ problems was 23 times shorter. For the tested problem with  $n=3038,~p=500$, run time without the short-cut would roughly be 500 times longer (and therefore not tested). Let $d_{ij}$ be the distance between demand points $i$ and $j$. We describe the short-cut for Algorithm~3b.

Three vectors of length $n$ each are saved in memory: $D^{(1)}$, $D^{(2)}$, and $D^{(3)}$.
\begin{enumerate}
\item Once the initial set of selected demand points $P$ is established (either by RAND or by CONS), the vector $D^{(1)} =\left\{ D^{(1)}_j,~j=1,\ldots,n\right\}$ is constructed giving the list of shortest distances between  each demand point $j$ and the $p$ facilities $\in P$.
\item \label{stp2a}Once the facility located at demand point $i\in P$ is removed from $P$, go over all demand points $j=1,\ldots,n$, and
\begin{enumerate}
\item If $d_{ij}>D^{(1)}_j$, set $D^{(2)}_j=D^{(1)}_j$.
\item Otherwise $(d_{ij}=D^{(1)}_j)$, find the minimum distance between demand point $j$ and the remaining facilities in $P\setminus \{i\}$ and assign it to  $D^{(2)}_j$.
\end{enumerate} 
\item\label{ST3} Suppose that demand point $k\in J\setminus P$ replaces demand point $i$. Go over all demand points $j=1,\ldots,n$, and 
\begin{enumerate}
\item If $d_{kj}\ge D^{(2)}_j$, set $D^{(3)}_j=D^{(2)}_j$.
\item If $d_{kj}< D^{(2)}_j$, set $D^{(3)}_j=d_{kj}$.
\end{enumerate}
\item The objective function is calculated using the vector $D^{(3)}$.
\item \label{ST5} If a better value of the objective function is found, set  $D^{(1)}=D^{(3)}$, replace demand point $i$ with demand point $k$, and return to Step \ref{stp2a} with the updated set $P$.
\item Continue the evaluation of all exchanges and if all fail to improve the objective function, stop.

\end{enumerate}

\subsection*{Notes}

\begin{enumerate}
	\item  In Step \ref{ST3}: Demand points in the remaining selected facilities in $P$ have a shortest distance of zero and can be skipped in the calculations but we did not observe a noticeable change in run time by implementing it.
\item  In Algorithm 3a continue the evaluations of the objective function in Step \ref{ST5} until all exchanges are evaluated and select the best improving pair $i,k$ to update $P$ and $D^{(1)}$.
\end{enumerate}

\begin{table}[ht!]
	\begin{center}
		\caption{\label{opt}Best known objectives of the test problems}
		\medskip
		\begin{tabular}{|c|c|r|c||c|c|r|c||c|c|r|c|}
			\hline
			$n$&$p$&objective&$\dagger$&$n$&$p$&objective&$\dagger$&$n$&$p$&objective&$\dagger$\\
			\hline
100&	5&	164.6011&	0.02&	400&	15&	362.7120&	0.39&	700&	25&	482.5661&	1.25\\
100&	10&	100.7650&	0.05&	400&	20&	304.1061&	0.54&	800&	5&	1372.8710&	0.23\\
100&	15&	74.4746&	0.09&	400&	25&	266.3945&	0.59&	800&	10&	928.7004&	0.59\\
100&	20&	59.4779&	0.21&	500&	5&	856.1153&	0.13&	800&	15&	743.1017&	0.93\\
100&	25&	49.1846&	0.37&	500&	10&	575.6737&	0.31&	800&	20&	633.9782&	1.36\\
200&	5&	329.0968&	0.04&	500&	15&	449.8948&	0.48&	800&	25&	557.1867&	1.46\\
200&	10&	213.1025&	0.09&	500&	20&	382.6915&	0.71&	900&	5&	1545.5993&	0.29\\
200&	15&	167.1654&	0.16&	500&	25&	337.3002&	0.80&	900&	10&	1053.7279&	0.69\\
200&	20&	140.0728&	0.23&	600&	5&	1030.9282&	0.17&	900&	15&	844.0657&	1.14\\
200&	25&	120.5562&	0.28&	600&	10&	694.2726&	0.41&	900&	20&	718.9711&	1.53\\
300&	5&	505.9990&	0.06&	600&	15&	547.8102&	0.61&	900&	25&	634.8785&	1.72\\
300&	10&	331.5499&	0.14&	600&	20&	460.6433&	0.87&	1000&	5&	1731.6308&	0.35\\
300&	15&	259.6754&	0.24&	600&	25&	408.3926&	0.95&	1000&	10&	1177.9664&	0.81\\
300&	20&	216.8050&	0.35&	700&	5&	1198.9113&	0.21&	1000&	15&	942.4672&	1.28\\
300&	25&	191.5259&	0.40&	700&	10&	807.4504&	0.47&	1000&	20&	798.5461&	1.73\\
400&	5&	685.1978&	0.10&	700&	15&	647.6007&	0.75&	1000&	25&	705.8626&	1.94\\
400&	10&	458.8549&	0.24&	700&	20&	548.0676&	1.12&	&	&	&	\\
			\hline							
\multicolumn{12}{l}{$\dagger$ Time in minutes per run.}
\end{tabular}
	\end{center}
\end{table}

\begin{table}[ht!]
	\begin{center}
		\caption{\label{Comp}Percent above best known solution of the best value obtained in 100 runs}
		\medskip
		\begin{tabular}{|c|c||c|c|c||c|c||c|c|c|}
			\hline
			$n$&$p$&ALT&IALT&RATIO&$n$&$p$&ALT&IALT&RATIO\\
			\hline
			100&	5&	0.01\%&	0.00\%&	0.00\%&	600&	5&	0.00\%&	0.00\%&	0.00\%\\
			100&	10&	0.80\%&	0.00\%&	0.17\%&	600&	10&	0.00\%&	0.00\%&	0.00\%\\
			100&	15&	2.48\%&	0.70\%&	0.00\%&	600&	15&	0.06\%&	0.13\%&	0.01\%\\
			100&	20&	1.60\%&	0.62\%&	0.62\%&	600&	20&	0.13\%&	0.01\%&	0.12\%\\
			100&	25&	4.27\%&	2.89\%&	1.10\%&	600&	25&	1.21\%&	0.68\%&	0.58\%\\
			200&	5&	0.00\%&	0.00\%&	0.00\%&	700&	5&	0.00\%&	0.00\%&	0.00\%\\
			200&	10&	0.00\%&	0.00\%&	0.00\%&	700&	10&	0.01\%&	0.00\%&	0.00\%\\
			200&	15&	0.10\%&	0.00\%&	0.57\%&	700&	15&	0.01\%&	0.00\%&	0.05\%\\
			200&	20&	1.52\%&	1.58\%&	1.17\%&	700&	20&	0.43\%&	0.22\%&	0.25\%\\
			200&	25&	4.21\%&	1.85\%&	2.10\%&	700&	25&	1.28\%&	0.47\%&	0.22\%\\
			300&	5&	0.00\%&	0.00\%&	0.00\%&	800&	5&	0.01\%&	0.00\%&	0.00\%\\
			300&	10&	0.01\%&	0.00\%&	0.00\%&	800&	10&	0.01\%&	0.00\%&	0.00\%\\
			300&	15&	0.86\%&	0.73\%&	0.73\%&	800&	15&	0.01\%&	0.01\%&	0.01\%\\
			300&	20&	1.82\%&	2.05\%&	0.01\%&	800&	20&	0.13\%&	0.07\%&	0.10\%\\
			300&	25&	1.95\%&	1.57\%&	0.96\%&	800&	25&	0.92\%&	0.55\%&	0.44\%\\
			400&	5&	0.00\%&	0.00\%&	0.00\%&	900&	5&	0.00\%&	0.00\%&	0.00\%\\
			400&	10&	0.00\%&	0.00\%&	0.00\%&	900&	10&	0.04\%&	0.00\%&	0.00\%\\
			400&	15&	0.02\%&	0.41\%&	0.09\%&	900&	15&	0.02\%&	0.00\%&	0.00\%\\
			400&	20&	0.24\%&	0.11\%&	0.06\%&	900&	20&	0.01\%&	0.30\%&	0.30\%\\
			400&	25&	1.07\%&	1.20\%&	0.03\%&	900&	25&	0.83\%&	0.59\%&	0.45\%\\
			500&	5&	0.00\%&	0.00\%&	0.00\%&	1000&	5&	0.00\%&	0.00\%&	0.00\%\\
			500&	10&	0.02\%&	0.00\%&	0.00\%&	1000&	10&	0.00\%&	0.00\%&	0.01\%\\
			500&	15&	0.21\%&	0.01\%&	0.00\%&	1000&	15&	0.02\%&	0.00\%&	0.03\%\\
			500&	20&	0.76\%&	0.06\%&	0.06\%&	1000&	20&	0.06\%&	0.02\%&	0.02\%\\
			500&	25&	0.89\%&	1.13\%&	0.50\%&	1000&	25&	0.42\%&	0.15\%&	0.02\%\\
			\hline							
		\end{tabular}
	\end{center}
\end{table}

\begin{table}[ht!]
	\begin{center}
		\caption{\label{Taly}Comparing RATIO to IALT for $n=3,038$ instances}
		\medskip
		\setlength{\tabcolsep}{5pt}
		\begin{tabular}{|c|c||c|c|c|c|c||c|c|c|c|}
			\hline
&Best&\multicolumn{5}{|c||}{IALT used in \cite{DTD18}}&\multicolumn{4}{|c|}{Using RATIO (same program)}\\
\cline{3-11}
$p$&Known$\dagger$&(1)&(2)&(3)&(4)&(5)&(1)&(2)&(3)&(4)\\
			\hline
50&	505,875.76&	10&	0\%&	0\%&	0.62&	0.35&	10&	0\%&	0\%&	0.29\\
100&	351,171.15&	9&	0\%&	0.002\%&	1.29&	0.72&	9&	0\%&	0.002\%&	0.66\\
150&	279,724.73&	8&	0\%&	0.003\%&	2.98&	1.66&	6&	0\%&	0.008\%&	1.33\\
200&	236,209.47&	0&	0.000\%&	0.007\%&	4.29&	2.39&	4&	0\%&	0.002\%&	2.39\\
250&	206,454.64&	1&	0\%&	0.002\%&	6.95&	3.86&	5&	0\%&	0.001\%&	3.82\\
300&	{\bf 184,799.90}&	0&	0.001\%&	0.008\%&	10.37&	5.76&	1&	0\%&	0.007\%&	5.93\\
350&	168,246.96&	1&	0\%&	0.008\%&	18.59&	10.33&	1&	0\%&	0.005\%&	12.84\\
400&	154,554.55&	0&	0.002\%&	0.011\%&	28.10&	15.61&	0&	0.000\%&	0.014\%&	15.34\\
450&	143,267.54&	0&	0.003\%&	0.013\%&	41.95&	23.30&	0&	0.000\%&	0.011\%&	20.94\\
500&	{\bf 133,547.50}&	0&	0.006\%&	0.017\%&	49.16&	27.31&	1&	0\%&	0.017\%&	35.33\\
\hline										
\multicolumn{2}{|l||}{Average:}&		2.9&	0.0012\%&	0.0072\%&	16.43&	9.13&	3.7&	0.0000\%&	0.0068\%&	9.89\\
			\hline
\multicolumn{11}{l}{$\dagger $ New best known solution marked in boldface.}\\							
\multicolumn{11}{l}{(1) Number of times in 10 runs that best known solution obtained.}\\							
\multicolumn{11}{l}{(2) Percent of best found solution above best known solution.}\\							
\multicolumn{11}{l}{(3) Percent of average solution above best known solution.}\\							
\multicolumn{11}{l}{(4) Run time in hours for one run.}\\							
\multicolumn{11}{l}{(5) Adjusted run time on the faster computer.}\\							
		\end{tabular}
	\end{center}
\end{table}

\begin{table}[ht!]
	\begin{center}
		\caption{\label{opt1}Comparing Average results of RATIO from Various Starting Solutions}
		\medskip
		\begin{tabular}{|c|c|c|c|c|c||c|c|c|c|c|c|}
			\hline
			$n$&$p$&(1)&(2)&(3)&(4)&$n$&$p$&(1)&(2)&(3)&(4)\\
			\hline
100&	5&	0.78\%&	0.93\%&	0.57\%&	0.88\%&	600&	5&	0.39\%&	0.53\%&	0.34\%&	0.27\%\\
100&	10&	4.57\%&	1.37\%&	0.26\%&	0.32\%&	600&	10&	1.47\%&	0.75\%&	0.47\%&	0.43\%\\
100&	15&	7.90\%&	3.62\%&	0.70\%&	0.73\%&	600&	15&	3.07\%&	1.70\%&	0.62\%&	0.55\%\\
100&	20&	10.23\%&	6.49\%&	0.33\%&	0.44\%&	600&	20&	4.20\%&	1.90\%&	0.54\%&	0.57\%\\
100&	25&	12.69\%&	4.12\%&	0.20\%&	0.19\%&	600&	25&	4.04\%&	2.33\%&	0.59\%&	0.60\%\\
200&	5&	1.08\%&	1.59\%&	0.24\%&	0.28\%&	700&	5&	0.62\%&	0.59\%&	0.35\%&	0.38\%\\
200&	10&	6.24\%&	1.51\%&	0.18\%&	0.11\%&	700&	10&	1.13\%&	0.55\%&	0.35\%&	0.31\%\\
200&	15&	6.70\%&	3.89\%&	0.54\%&	0.56\%&	700&	15&	2.20\%&	1.24\%&	0.42\%&	0.56\%\\
200&	20&	7.10\%&	4.48\%&	0.36\%&	0.37\%&	700&	20&	3.44\%&	1.44\%&	0.53\%&	0.57\%\\
200&	25&	8.63\%&	4.19\%&	0.63\%&	0.46\%&	700&	25&	3.80\%&	2.16\%&	0.74\%&	0.80\%\\
300&	5&	1.09\%&	1.57\%&	0.26\%&	0.25\%&	800&	5&	0.41\%&	0.50\%&	0.43\%&	0.33\%\\
300&	10&	4.22\%&	0.97\%&	0.48\%&	0.48\%&	800&	10&	0.66\%&	0.44\%&	0.21\%&	0.16\%\\
300&	15&	6.36\%&	3.54\%&	0.67\%&	0.66\%&	800&	15&	2.64\%&	1.38\%&	0.56\%&	0.56\%\\
300&	20&	7.15\%&	4.10\%&	0.27\%&	0.33\%&	800&	20&	2.56\%&	1.24\%&	0.34\%&	0.34\%\\
300&	25&	7.57\%&	3.28\%&	0.43\%&	0.38\%&	800&	25&	3.47\%&	2.14\%&	0.60\%&	0.65\%\\
400&	5&	1.16\%&	1.24\%&	0.47\%&	0.53\%&	900&	5&	0.58\%&	0.61\%&	0.41\%&	0.41\%\\
400&	10&	1.70\%&	0.76\%&	0.52\%&	0.43\%&	900&	10&	0.42\%&	0.46\%&	0.30\%&	0.32\%\\
400&	15&	3.40\%&	1.69\%&	0.41\%&	0.42\%&	900&	15&	2.38\%&	1.30\%&	0.71\%&	0.74\%\\
400&	20&	4.80\%&	2.59\%&	0.26\%&	0.32\%&	900&	20&	2.83\%&	1.53\%&	0.69\%&	0.62\%\\
400&	25&	5.81\%&	2.90\%&	0.45\%&	0.42\%&	900&	25&	3.27\%&	1.88\%&	0.69\%&	0.68\%\\
500&	5&	0.71\%&	0.74\%&	0.29\%&	0.33\%&	1000&	5&	0.67\%&	0.80\%&	0.52\%&	0.56\%\\
500&	10&	1.39\%&	0.80\%&	0.44\%&	0.53\%&	1000&	10&	0.86\%&	0.67\%&	0.42\%&	0.53\%\\
500&	15&	4.35\%&	1.62\%&	0.40\%&	0.34\%&	1000&	15&	1.44\%&	0.96\%&	0.40\%&	0.40\%\\
500&	20&	4.68\%&	1.72\%&	0.37\%&	0.41\%&	1000&	20&	2.37\%&	1.65\%&	0.45\%&	0.43\%\\
500&	25&	4.91\%&	2.61\%&	0.55\%&	0.52\%&	1000&	25&	2.57\%&	1.63\%&	0.49\%&	0.51\%\\
			\hline							
			\multicolumn{12}{l}{(1) Average of 100 RAND solutions above best known solution.}\\
			\multicolumn{12}{l}{(2) Average of 100 CONS solutions above best known solution.}\\
		\multicolumn{12}{l}{(3) Average of 100 DESC solutions above best known solution.}\\
		\multicolumn{12}{l}{(4) Average of 100 COMB solutions above best known solution.}
			
		\end{tabular}
	\end{center}
\end{table}

\begin{table}[ht!]
	\begin{center}
		\caption{\label{Best}Comparing Best Results of Various Starting Solutions for Small Problems}
		\medskip
		\begin{tabular}{|c|c||c|c||c|c||c|c||c|c|}
			\hline
			&&\multicolumn{2}{|c||}{RAND}&\multicolumn{2}{|c||}{CONS}&\multicolumn{2}{|c||}{DESC}&\multicolumn{2}{|c|}{COMB}\\
			\cline{3-10}
			$n$&$p$&(1)&(2)&(1)&(2)&(1)&(2)&(1)&(2)\\
			\hline
100&	5&	0\%&	7&	0\%&	13&	0\%&	8&	0\%&	2\\
100&	10&	0\%&	3&	0\%&	10&	0\%&	33&	0\%&	16\\
100&	15&	0.64\%&	0&	0.44\%&	0&	0\%&	14&	0\%&	8\\
100&	20&	1.17\%&	0&	1.91\%&	0&	0\%&	68&	0\%&	53\\
100&	25&	1.70\%&	0&	0.50\%&	0&	0\%&	34&	0\%&	36\\
200&	5&	0\%&	34&	0\%&	33&	0\%&	38&	0\%&	48\\
200&	10&	0\%&	6&	0\%&	5&	0\%&	24&	0\%&	29\\
200&	15&	0.10\%&	0&	0\%&	1&	0\%&	9&	0\%&	7\\
200&	20&	1.29\%&	0&	0.34\%&	0&	0.04\%&	0&	0\%&	2\\
200&	25&	1.33\%&	0&	0.90\%&	0&	0.01\%&	0&	0.01\%&	0\\
300&	5&	0\%&	4&	0\%&	12&	0\%&	24&	0\%&	34\\
300&	10&	0\%&	11&	0\%&	22&	0\%&	41&	0\%&	46\\
300&	15&	0.92\%&	0&	0.00\%&	0&	0\%&	27&	0\%&	33\\
300&	20&	0.01\%&	0&	0.58\%&	0&	0\%&	11&	0\%&	17\\
300&	25&	1.47\%&	0&	0.19\%&	0&	0\%&	3&	0\%&	3\\
400&	5&	0\%&	33&	0\%&	28&	0\%&	47&	0\%&	42\\
400&	10&	0\%&	1&	0\%&	2&	0.00\%&	0&	0.01\%&	0\\
400&	15&	0.10\%&	0&	0.00\%&	0&	0.00\%&	0&	0.00\%&	0\\
400&	20&	0.23\%&	0&	0.03\%&	0&	0.02\%&	0&	0.02\%&	0\\
400&	25&	0.40\%&	0&	0.30\%&	0&	0\%&	3&	0\%&	1\\
500&	5&	0\%&	40&	0\%&	35&	0\%&	72&	0\%&	67\\
500&	10&	0\%&	2&	0\%&	2&	0\%&	1&	0.00\%&	0\\
500&	15&	0.02\%&	0&	0.00\%&	0&	0\%&	36&	0\%&	28\\
500&	20&	0.46\%&	0&	0.34\%&	0&	0.02\%&	0&	0.02\%&	0\\
500&	25&	1.53\%&	0&	1.02\%&	0&	0\%&	2&	0\%&	1\\
\hline									
\multicolumn{2}{|l||}{Average:}&		0.455\%&	5.64&	0.263\%&	6.52&	0.004\%&	19.80&	0.003\%&	18.92\\
			\hline							
			\multicolumn{10}{l}{(1) Percent of best found solution above best known solution.}\\
			\multicolumn{10}{l}{(2) Number of times best known solution found.}
		\end{tabular}
	\end{center}
\end{table}

\begin{table}[ht!]
	\begin{center}
		\caption{\label{BestL}Comparing Best Results of Various Starting Solutions for Large Problems}
		\medskip
		\begin{tabular}{|c|c||c|c||c|c||c|c||c|c|}
			\hline
			&&\multicolumn{2}{|c||}{RAND}&\multicolumn{2}{|c||}{CONS}&\multicolumn{2}{|c||}{DESC}&\multicolumn{2}{|c|}{COMB}\\
			\cline{3-10}
			$n$&$p$&(1)&(2)&(1)&(2)&(1)&(2)&(1)&(2)\\
			\hline
600&	5&	0\%&	5&	0.00\%&	0&	0\%&	6&	0\%&	15\\
600&	10&	0\%&	2&	0\%&	8&	0\%&	9&	0\%&	11\\
600&	15&	0.01\%&	0&	0.03\%&	0&	0.01\%&	0&	0.01\%&	0\\
600&	20&	0.01\%&	0&	0.01\%&	0&	0.00\%&	0&	0\%&	1\\
600&	25&	0.99\%&	0&	0.66\%&	0&	0\%&	1&	0\%&	1\\
700&	5&	0\%&	27&	0\%&	27&	0\%&	47&	0\%&	48\\
700&	10&	0\%&	15&	0\%&	12&	0\%&	20&	0\%&	19\\
700&	15&	0.01\%&	0&	0.01\%&	0&	0\%&	6&	0\%&	5\\
700&	20&	0.49\%&	0&	0.08\%&	0&	0.03\%&	0&	0.03\%&	0\\
700&	25&	0.26\%&	0&	0.42\%&	0&	0.01\%&	0&	0\%&	1\\
800&	5&	0.00\%&	0&	0.00\%&	0&	0.00\%&	0&	0.00\%&	0\\
800&	10&	0\%&	5&	0\%&	2&	0\%&	1&	0.00\%&	0\\
800&	15&	0.02\%&	0&	0.01\%&	0&	0.01\%&	0&	0.01\%&	0\\
800&	20&	0.02\%&	0&	0.19\%&	0&	0.05\%&	0&	0.02\%&	0\\
800&	25&	0.68\%&	0&	0.20\%&	0&	0\%&	1&	0.01\%&	0\\
900&	5&	0\%&	8&	0\%&	15&	0\%&	13&	0\%&	21\\
900&	10&	0\%&	1&	0\%&	4&	0\%&	5&	0\%&	3\\
900&	15&	0\%&	1&	0\%&	2&	0\%&	24&	0\%&	26\\
900&	20&	0.23\%&	0&	0\%&	1&	0.00\%&	0&	0.00\%&	0\\
900&	25&	0.12\%&	0&	0.42\%&	0&	0\%&	2&	0\%&	2\\
1000&	5&	0\%&	5&	0\%&	6&	0\%&	17&	0\%&	15\\
1000&	10&	0\%&	1&	0\%&	1&	0.00\%&	0&	0.00\%&	0\\
1000&	15&	0.03\%&	0&	0.00\%&	0&	0\%&	2&	0\%&	1\\
1000&	20&	0\%&	1&	0.36\%&	0&	0\%&	1&	0\%&	1\\
1000&	25&	0.15\%&	0&	0.07\%&	0&	0.01\%&	0&	0.02\%&	0\\
\hline									
\multicolumn{2}{|l||}{Average:}&		0.121\%&	2.84&	0.098\%&	3.12&	0.005\%&	6.20&	0.004\%&	6.80\\
\hline									
\multicolumn{2}{|l||}{Total Average:}&		0.288\%&	4.24&	0.181\%&	4.82&	0.004\%&	13.00&	0.003\%&	12.86\\
			\hline							
			\multicolumn{10}{l}{(1) Percent of best found solution above best known solution.}\\
			\multicolumn{10}{l}{(2) Number of times best known solution found.}
		\end{tabular}
	\end{center}
\end{table}

\begin{table}[ht!]
	\begin{center}
		\caption{\label{summary}Comparing performance of the starting solutions methods}
		\medskip
\setlength{\tabcolsep}{5pt}
		\begin{tabular}{|l||c|c|c|c|}
			\hline
&	RAND&	CONS&	DESC&	COMB\\
Property&	1000 Runs&	1000 Runs&	100 Runs&	100 Runs\\
\hline				
Average of best found above best known&	0.088\%&	0.039\%&	0.004\%&	0.003\%\\
Average of average above best known&	3.53\%&	1.85\%&	0.45\%&	0.46\%\\
Total run time (sec) for $n=100,~ p=5$&	0.19&	0.17&	0.06&	0.05\\
Total run time (sec) for $n=1000,~ p=25$&	15.34&	13.58&	22.62&	23.41\\
Average total run time (sec) for all runs&	3.76&	3.30&	4.80&	4.87\\
\hline				
		\end{tabular}
	\end{center}
\end{table}

\section{Computational Experiments}

Computer programs were coded in Fortran using double precision arithmetic. The programs were compiled by an Intel 11.1 Fortran Compiler with no parallel processing. They were run on a desktop with the Intel i7-6700 3.4GHz CPU processor and 16GB RAM. Only one thread was used. For comparison, the results reported in \citet{DTD18} were run on a slower desktop (also only one thread used) with the Intel 870/i7 2.93GHz CPU Quad processor and 8GB RAM. Run times on the slower computer were longer by about 80\%.

To allow for easy replication for future comparisons, we tested instances that were randomly generated by the method proposed in \cite{DKS18,DDK19} (detailed in the Appendix), with all weights set to 1. 1,000 points were generated, and for $n<1,000$ the first $n$ points were used.
For reference purposes we depict in Table \ref{opt} the best found solutions by the method in \cite{DTD18} which is based on the algorithm in \cite{DS15,DBMS14} using a specially designed population of starting solutions, followed by a genetic algorithm \citep{H75,G06} followed by a variable neighborhood search \citep{HM97, MH97}. Each instance was solved 10 times and in all cases all 10 replications yielded the same solution. Total time for all 50 instances run 10 times each was less than 5 hours.

\subsection{Comparing IALT to RATIO}

We first compared the original ALT algorithm to IALT, and its RATIO version. Starting solutions were randomly generated by the same starting seed. However, since every algorithm generated a different number of random numbers, starting solutions of the second case and higher are different for each algorithm. In Table~\ref{Comp} we show the percent above the best known solution of the best result by ALT, IALT, and RATIO. The RATIO variant performed significantly better, especially for larger values of $p$. In nine of the ten $p=25$ instances RATIO was better than IALT. A paired t-test comparing IALT to RATIO for all 50 instances has a p-value of 0.012. The averages are significantly better as well with p-value of 0.028.  Both IALT and RATIO are clearly better than ALT.
Run times for IALT or RATIO are very fast. One run of the largest problem required about 0.015 seconds.

To further explore whether RATIO outperforms IALT, we examined a challenging problem set with $n=3,038$  \citep{Re91}  and $p=50,100,\ldots,500$, that has also been tested in previous papers. The results are reported in Table \ref{Taly}. RATIO indeed outperformed IALT, and was even able to find two new best known solutions. The best found solution for the two instances where it did not find the best known solution, $p=400,450$, were only 0.0002\% and 0.0003\% above the best known solution. The best found solutions are significantly better by RATIO with p-value of 0.041. The average solutions are also better but not within statistical significance (p-value=0.32). 
We believe that RATIO accommodates outliers, which have a larger distance to the closest facility,  better than the original IALT, because such outliers are more likely to have the $L$ smallest ratios and less likely to have the $L$ smallest differences.

\subsection{Comparing Various Starting Solution Approaches}

Four heuristics derived from four starting solutions followed by RATIO, are compared: (i) random selection of $p$ demand points (RAND), (ii) the construction method for selecting $p$ demand points (CONS), (ii) the descent algorithm applied on the random start (DESC), and (iv) the combined approach (COMB) which is DESC applied on CONS. There are two variants of DESC labeled Algorithm 3a and Algorithm 3b. Algorithm 3b performed better in both quality of the solution and run time. Therefore, we do not report any results of Algorithm 3a, that is, DESC and COMB both used Algorithm~3b.

The percentage of the average results above the best known solution are depicted in Table \ref{opt1}. The CONS algorithm clearly outperformed the RAND algorithm. The descent approach (DESC or COMB) clearly improves the average solution further but DESC and COMB performed about equally well.

In Tables \ref{Best} and \ref{BestL} the percentage of the best found solution above the best known solution, and the number of times the best known solution was found in 100 runs, are reported. As should be expected, DESC and COMB clearly outperformed RAND and CONS. However run times are longer for DESC and COMB.

Since run times for DESC and COMB are up to 10 times longer, for fair comparison of the best found solutions we compare the best results of 1000 runs of RAND and CONS with only 100 runs of the DESC and the COMB.  A short summary of the results is depicted in Table \ref{summary}. DESC and COMB still performed much better under these conditions. Note that without the short-cut proposed in Section \ref{short}, run times of DESC and COMB are about 23 times longer and to achieve similar run times for these problems RAND and COMB should have been run 23,000 times for a fair comparison.

\begin{table}[ht!]
	\begin{center}
		\caption{\label{3038C}Solving the  $n=3,038$ Instances using RATIO}
		\medskip
		\setlength{\tabcolsep}{2.5pt}
		\begin{tabular}{|c|c||c|c|c||c|c|c||c|c|c||c|c|c|}
			\hline
			&Best&\multicolumn{3}{|c||}{RAND 100 Runs}&\multicolumn{3}{|c||}{CONS 100 Runs}&\multicolumn{3}{|c||}{DESC 10 Runs}&\multicolumn{3}{|c|}{COMB 10 Runs}\\
			\cline{3-14}
			$p$&Known&(1)&(2)&(3)&(1)&(2)&(3)&(1)&(2)&(3)&(1)&(2)&(3)\\
			\hline
50&	505,875.76&	0.64\%&	2.59\%&	0.34&	0.94\%&	1.86\%&	0.26&	0.32\%&	0.55\%&	0.98&	0.20\%&	0.59\%&	0.81\\
100&	351,171.15&	1.87\%&	3.64\%&	0.97&	0.92\%&	2.46\%&	0.83&	0.60\%&	0.84\%&	1.78&	0.58\%&	0.81\%&	1.63\\
150&	279,724.73&	2.51\%&	4.64\%&	1.68&	1.86\%&	2.92\%&	1.42&	0.33\%&	0.74\%&	2.51&	0.70\%&	0.86\%&	2.92\\
200&	236,209.47&	3.48\%&	5.73\%&	2.50&	2.00\%&	3.29\%&	2.02&	0.57\%&	0.70\%&	3.43&	0.59\%&	0.83\%&	3.46\\
250&	206,454.64&	3.87\%&	6.16\%&	3.31&	2.38\%&	3.75\%&	2.66&	0.62\%&	0.74\%&	4.31&	0.50\%&	0.73\%&	4.05\\
300&	184,799.90&	5.00\%&	6.87\%&	4.18&	3.36\%&	4.05\%&	3.29&	0.63\%&	0.85\%&	4.96&	0.76\%&	0.89\%&	4.87\\
350&	168,246.96&	5.71\%&	7.48\%&	4.89&	3.32\%&	4.16\%&	3.89&	0.68\%&	0.84\%&	5.53&	0.61\%&	0.85\%&	5.75\\
400&	154,554.55&	6.37\%&	8.18\%&	5.67&	3.57\%&	4.51\%&	4.52&	0.83\%&	0.98\%&	6.52&	0.74\%&	0.97\%&	6.39\\
450&	143,267.54&	6.51\%&	8.63\%&	6.44&	3.88\%&	4.71\%&	5.00&	0.75\%&	0.91\%&	7.36&	0.65\%&	0.88\%&	7.31\\
500&	133,547.50&	7.87\%&	9.69\%&	7.17&	4.20\%&	5.06\%&	5.50&	0.79\%&	0.96\%&	8.28&	0.83\%&	0.96\%&	7.72\\
\hline													
\multicolumn{2}{|l||}{Average:}&		4.38\%&	6.36\%&	3.71&	2.64\%&	3.68\%&	2.94&	0.61\%&	0.81\%&	4.57&	0.62\%&	0.84\%&	4.49\\
			\hline
			\multicolumn{14}{l}{(1) Percent of best found solution above best known solution.}\\							
			\multicolumn{14}{l}{(2) Percent of average solution above best known solution.}\\							
			\multicolumn{14}{l}{(3) Run time in minutes for all runs.}\\							
		\end{tabular}
	\end{center}
\end{table}

\subsection{Solving the $n=3,038$ Instances}

The methods suggested in \cite{DTD18,DS15,DBMS14} are designed to find good solutions but require a very long run time. They are based on complicated starting solutions followed by a genetic algorithm and then a variable neighborhood search. The $n=3,038,~p=500$ instance was solved in about 35 hours for one run. It took about seven weeks to get the results depicted in Table \ref{Taly}. Using the suggested starting solutions (DESC or COMB) and applying the RATIO variant of IALT on the generated starting solutions takes less than 0.8 minutes per run for this large problem, which is almost 3000 times faster. If one is interested in getting reasonably good results in a relatively short computational time, the procedures proposed in this paper can be used. Also, run time of COMB or DESC is about proportional to $p$ while run time for the methods suggested in \cite{DTD18,DS15,DBMS14} is proportional to about $p^{2.1}$.

We solved the $n=3,038$ instances by the four approaches and report in Table \ref{3038C} the best and average results obtained by solving each instance 100 times by RAND and CONS and 10 times by DESC and COMB. The best and even average results are below 1\% above the best known solution for all instances.  Run times for the largest problem are about 5-8 minutes for all runs. They are about 3 to 4 seconds for one run of the largest problem by RAND and CONS and about 50 seconds for one run of DESC and COMB.

\section{Conclusions}

It is well known that good solutions of the discrete $p$-median problem (where facility sites are restricted to the set of demand points) are also generally good solutions of the continuous counterpart, and hence, may be used as a component of heuristic algorithms for solving the planar $p$-median problem
\citep[e.g., see ][]{HMT98,BHMT00,BDMS13}.
The main approach used in previous work has been to first solve the discrete problem exactly (which can be very time-consuming), and then perform a continuous adjustment step on the obtained partition of the demand points to finish with a good solution to the original problem. Here we wish to further explore the use of “good” discrete starting solutions. Two new approaches for finding such solutions are proposed. The first uses a constructive method with random component referred to as CONS; the second applies a simple local search (vertex swap), referred to as DESC, to randomly generated selections of $p$ demand points (RAND). A third approach (COMB) combines CONS and DESC. These discrete solutions are then used as starting solutions in a powerful local search that iterates between Cooper's famous locate-allocate method and a transfer follow-up step \citep{BD12}. A new selection rule for the transfer follow-up is also tested yielding improved results.

Extensive computational experiments on medium and large-sized problem sets show that the CONS, DESC and COMB approaches significantly improve the quality of obtained solutions compared to starting with randomly-generated solutions (RAND). Furthermore, DESC and COMB perform substantially better than CONS. The heuristics developed here also perform well compared to state-of-the-art meta heuristics from the literature. For the largest problem set investigated ($n = 3038$ demand points), the average results of DESC and COMB for all instances tested were less than 1 \% above the best-known solutions, while using a miniscule fraction of the computing time of these sophisticated algorithms.

\begin{center}{\large\bf Appendix: Generating Random Configurations}\end{center}

We generate a sequence of integer numbers in the open range (0, 100,000). A starting seed $r_1$, which is the first number in the sequence, is selected. The sequence is generated by the following rule for $k\ge 1$:
\begin{itemize}\narrow
	\item Set $\theta=12219 r_k$.
	\item Set $r_{k+1}=\theta-\lfloor\frac{\theta}{100000}\rfloor\times 100000$, i.e., $r_{k+1}$ is the remainder of dividing $\theta$ by 100000.
\end{itemize}

For the $x$ coordinates we used $r_1=97$ and for the $y$-coordinates we used $r_1=367$.  To define the coordinates we divide $r_k$ by 10000.
 
		\renewcommand{\baselinestretch}{1}
		\renewcommand{\arraystretch}{1}
		\large
		\normalsize

		\bibliographystyle{apalike}
		

\end{document}